\documentclass[a4paper,11pt]{article}

\usepackage{amsmath}
\usepackage{amssymb}
\usepackage{graphicx}

\newtheorem{theorem}{Theorem}[section]
\newtheorem{lemma}{Lemma}[section]
\newtheorem{proposition}{Proposition}[section]
\newtheorem{corollary}{Corollary}[section]
\newtheorem{example}{Example}[section]

\begin{document}
\title{Diagonal Riccati Stability and Applications}

\author{Alexander Aleksandrov\thanks{Faculty of Applied Mathematics and Control Processes, St. Petersburg State University, 35 Universitetskij Pr., 198504 Petrodvorets, St. Petersburg, Russia. email: alex43102006@yandex.ru} and Oliver Mason\thanks{Corresponding Author. Dept. of Mathematics and Statistics/Hamilton Institute, Maynooth University-National University of Ireland Maynooth, Maynooth, Co. Kildare, Ireland. email: oliver.mason@nuim.ie}}

\maketitle

\begin{abstract}
We consider the question of diagonal Riccati stability for a pair of real matrices $A, B$.  A necessary and sufficient condition for diagonal Riccati stability is derived and applications of this to two distinct cases are presented.  We also describe some motivations for this question arising in the theory of generalised Lotka-Volterra systems.
\end{abstract}

\begin{center}
\end{center}

\section{Introduction and Preliminaries}
We consider the following problem.  Given $A, B \in \mathbb{R}^{n \times n}$, determine conditions for the existence of \emph{diagonal} positive definite matrices $P, Q$ in $\mathbb{R}^{n\times n}$ such that 
\begin{equation}
\label{eq:Ricc3}
A^TP+PA + Q + PBQ^{-1}B^TP \prec 0,
\end{equation}
where $M \prec 0$ denotes that $M=M^T$ is negative definite.  Throughout the paper, $M \succ 0$ ($M \succeq 0$) denotes that $M$ is positive definite (positive semi-definite); $M \preceq 0$ denotes that $M$ is negative semi-definite.

When diagonal positive definite solutions $P, Q$ of \eqref{eq:Ricc3} exist, we say that the pair $A, B$ is \emph{diagonally Riccati stable}. 

Our interest in the question stems from the stability of time-delay systems.  
Specifically, when such a pair $P, Q$ exists, the linear 
time-delay system associated 
with $A, B$ admits a Lyapunov-Krasovskii functional of a particularly simple 
form  \cite{HALE}.  The more general question of when positive definite solutions of 
\eqref{eq:Ricc3}, not necessarily diagonal, exist  was highlighted 
in \cite{VERR} and some preliminary results on this question were also described in this reference.

We now introduce some notation and terminology, and recall some basic results that will be needed later in the paper.  

A matrix $A$ in $\mathbb{R}^{n \times n}$ is \emph{Metzler} if its off-diagonal elements are nonnegative: formally, $a_{ij} \geq 0$ for $i \neq j$.  For vectors $v, w$ in $\mathbb{R}^n$, $v \geq w$ is understood componentwise and means that $v_i \geq w_i$ for $1 \leq i\leq n$.  Similarly, $v > w$ means $v \geq w$, $v \neq w$ and $v \gg w$ means $v_i > w_i$ for $1 \leq i \leq n$.  

For a positive integer $n$, we denote by $Sym(n, \mathbb{R})$ the space of $n \times n$ real symmetric matrices.  For $A \in \mathbb{R}^{n \times n}$, $A^T$ denotes the tranpose of $A$.  A matrix $A \in \mathbb{R}^{n \times n}$ is \emph{Hurwitz} if all of its eigenvalues have negative real parts.  It is classical that $A$ is Hurwitz if and only if there exists some $P \succ 0$ with $A^TP+PA \prec 0$. 

The following result recalls some well-known facts concerning Metzler matrices \cite{HJ2}.
\begin{proposition}
\label{prop:Met1}
Let $A \in \mathbb{R}^{n \times n}$ be Metzler.  The following are equivalent:
\begin{itemize}
\item[(i)] $A$ is Hurwitz;
\item[(ii)] there exists some vector $v \gg 0$ in $\mathbb{R}^n$ with $Av \ll 0$;
\item[(iii)] there exists some positive definite diagonal matrix $D \in \mathbb{R}^{n \times n}$ with $A^TD+DA \prec 0$;
\item[(iv)] for every non-zero $v \in \mathbb{R}^{n}$, there is some index $i$ with $v_i(Av)_i < 0$.
\end{itemize}
\end{proposition}

The next result concerning the Riccati equation is based on the Schur complement and follows from Theorem 7.7.6 in \cite{HJ}. 
\begin{lemma}
\label{lem:Sch}
Let $A \in \mathbb{R}^{n \times n}$, $B \in \mathbb{R}^{n \times n}$ be given.  Then
the matrices $P\succ 0$, $Q\succ 0$ in $\mathbb{R}^{n \times n}$ satisfy (\ref{eq:Ricc3}) if and only if:
\begin{equation}
\label{eq:Schur}
S:= \left( \begin{array}{c c}
					A^TP+PA + Q & PB \\
					B^TP & -Q 
			\end{array} \right) \prec 0.
\end{equation}
\end{lemma}

\begin{theorem} \cite{BAR}
\label{thm:Sep} Let $C_1$, $C_2$ be non-empty convex subsets of a 
Euclidean Space $E$ with inner product $\langle \cdot, \cdot \rangle$.  
Further, assume that $C_2$ is a cone and $C_1 \cap C_2$ is empty. 
Then there exists a non-zero vector $v \in E$ such that 
\begin{eqnarray*}
\langle v, x \rangle &\geq & 0 \mbox{ for all } x \in C_1, \\
\langle v, x \rangle &\leq & 0 \mbox{ for all } x \in C_2. \\
\end{eqnarray*}
\end{theorem}

In this paper, we apply the above result to convex sets in the Euclidean space $Sym(m, \mathbb{R})$ of $m \times m$ real symmetric matrices equipped with the inner product $$\langle A, B \rangle = \textrm{trace}(AB).$$

\section{Diagonal Solutions to the Riccati Equation}
\label{sec:RiccDiag}
In a recent paper \cite{MAS}, it was shown that if $A$ is Metzler, $B$ 
is nonnegative and $A+B$ is Hurwitz, then there exist $P \succ 0$ diagonal 
and $Q \succ 0$ satisfying \eqref{eq:Ricc3} ($Q$ need not be diagonal).  Later, we shall show that it is in fact possible to choose $Q$ to be diagonal in this case.  

We first present necessary and sufficient conditions for the pair $A, B$ to be diagonally Riccati stable.  This result echoes the condition for the existence of diagonal solutions to the Lyapunov equation given in \cite{BBP}: see also \cite{SHO1, SHO2, KasBha} for more recent work on diagonal solutions of the Lyapunov equation.

\begin{theorem}
\label{thm:TM} Let $A, B \in \mathbb{R}^{n \times n}$ be given.  The following are equivalent.
\begin{itemize}
\item[(i)] There exist $P \succ 0$, $Q \succ 0$ diagonal satisfying \eqref{eq:Ricc3}.
\item[(ii)] For every non-zero positive semi-definite
\begin{equation}
\label{eq:H}
H = \left(\begin{array}{c c}
		H_{11} & H_{12} \\
		H_{12}^T & H_{22}
	  \end{array}\right)
\end{equation}
in $\mathbb{R}^{2n \times 2n}$ with $\textrm{diag}(H_{11}) \geq \textrm{diag}(H_{22})$, the matrix
$$AH_{11} + BH_{12}^T$$
has a negative diagonal entry.
\end{itemize}
\end{theorem}
\textbf{Proof:}
\emph{$(i) \Rightarrow (ii)$}:

If there exist diagonal matrices $P \succ 0, Q \succ 0$ satisfying \eqref{eq:Ricc3}, it is immediate from Lemma \ref{lem:Sch} that the matrix 
\begin{equation}
\label{eq:block}
R := \left(\begin{array}{c c}
		A^TP+PA + Q & PB \\
		B^TP & -Q
	\end{array}\right) \prec 0.
\end{equation}
Let a non-zero $H \in \mathbb{R}^{2n \times 2n}$ as in \eqref{eq:H} be given satisfying $\textrm{diag}(H_{11}) \geq \textrm{diag}(H_{22})$.  As $R$ is negative definite and $H$ is non-zero, it follows that 
\begin{equation}
\label{eq:TM11}
\textrm{trace}(HR) < 0.
\end{equation}
Expanding the left hand side of  \eqref{eq:TM11} we see that 
\begin{equation}
\label{eq:TM12}
\textrm{trace}(H_{11}A^TP + H_{11}PA + H_{11}Q + H_{12} B^TP + H_{12}^TPB - H_{22} Q) < 0
\end{equation}
which implies that
$$\textrm{trace}(2PAH_{11} + 2PBH_{12}^T + Q(H_{11} - H_{22}) ) < 0.$$
As $\textrm{diag}(H_{11}) \geq \textrm{diag}(H_{22})$ and $Q \succ 0$ is 
diagonal, it follows that $\textrm{trace}(Q(H_{11} - H_{22}) ) \geq 0$ and thus we can conclude that 
$$\textrm{trace}(P(AH_{11} + BH_{12}^T) ) < 0.$$
As $P \succ 0$ is diagonal, this implies that the matrix 
$$AH_{11} + BH_{12}^T$$ has at least one negative diagonal entry.  

\emph{$(ii) \Rightarrow (i)$}:  We shall prove the contrapositive.  Suppose that there exist no positive definite diagonal matrices $P$, $Q$ satisfying \eqref{eq:Ricc3}.  Consider the sets 
\begin{equation}
\label{eq:Cones1} \mathcal{C}_1 := \left \{ \left(\begin{array}{c c}
		A^TP+PA + Q & PB \\
		B^TP & -Q
	\end{array}\right) \mid Q, P \succ 0 \mbox{ diagonal } \in \mathbb{R}^{n \times n} \right\},
\end{equation}
and
$$\mathcal{C}_2 := \{H \in \mathbb{R}^{2n \times 2n} \mid H= H^T \prec 0\}.$$
Then it is easy to see that $\mathcal{C}_1$ and $\mathcal{C}_2$ are convex cones and by assumption we have that 
$$\mathcal{C}_1 \cap \mathcal{C}_2 = \emptyset.$$
It follows from Theorem \ref{thm:Sep} that there exists a non-zero  
matrix $H$ in $\textrm{Sym}(2n, \mathbb{R})$ such that
$$\textrm{trace}(H X) \leq 0 \mbox{ for all } X \in \mathcal{C}_2$$
$$\textrm{trace}(HX ) \geq 0 \mbox{ for all } X \in \mathcal{C}_1.$$
The first of the above two inequalities implies that $H \succeq 0$.  The second implies that 
\begin{equation}
\label{eq:sufftrace1} 
\textrm{trace}(H_{11}A^TP+ H_{11}PA + H_{12}B^TP + H_{12}^T PB + (H_{11} - H_{22}) Q ) \geq 0
\end{equation}
for all choices of diagonal matrices $P \succ 0$, $Q \succ 0$.  Rearranging, we see that for all $P \succ 0$, $Q \succ 0$ diagonal 
\begin{equation}
\label{eq:sufftrace2}
\textrm{trace}(2P(AH_{11} + BH_{12}^T)) \geq \textrm{trace}(Q (H_{22} - H_{11})).
\end{equation}
If there exists a diagonal $Q \succ 0$ such that 
$\textrm{trace}(Q(H_{22} - H_{11})) >0$,   then for a given diagonal $P \succ 0$, we could obtain a contradiction to \eqref{eq:sufftrace2} by considering $tQ$ for $t > 0$ sufficiently large.  This implies that $\textrm{trace}(Q(H_{22} - H_{11})) \leq 0$ for all diagonal $Q \succ 0$.  It follows from this that $\textrm{diag}(H_{22}) \leq \textrm{diag}(H_{11})$.  Moreover, for every diagonal $P \succ 0$, we must have 
$$\textrm{trace}(P(AH_{11} + BH_{12}^T)) \geq 0.$$
Now suppose that $AH_{11} + BH_{12}^T$ has a negative diagonal entry: say in the $i$th position.  Then we could choose $P \succ 0$ diagonal such that
$\textrm{trace}(P(AH_{11} + BH_{12}^T)) < 0 $ by setting the $i$th entry of $P$ to be large enough and keeping the other entries of $P$ small and positive.  It follows that under the assumption that the pair $A, B$ is not diagonally Riccati stable, that all diagonal entries of $AH_{11} + BH_{12}^T$ are nonnegative.  Hence, if $AH_{11} + BH_{12}^T$ has a negative diagonal entry for every $H$ satisfying the conditions of the theorem, there must exist diagonal $P \succ 0$, $Q \succ 0$ satisfying \eqref{eq:Ricc3}.  This completes the proof.

\vspace{5mm}

A simple application of Schur complements shows that if the pair $A, B$ is 
diagonally Riccati stable and $P$, $Q$ satisfy \eqref{eq:Ricc3}, then 
\begin{eqnarray*}
A^T P+P A &\prec& 0\\
(A+B)^TP +P(A+B) &\prec& 0.
\end{eqnarray*}

In the paper \cite{MAS}, it was asked whether or not the converse of this statement held.  It is possible to use Theorem \ref{thm:TM} to answer this and show that it is not the case.
\begin{example}
\label{ex:Ctr1} Consider 
$$A = \left(\begin{array}{c c}
		-1 & 0 \\
		-2 & -1
	\end{array}\right), \;\;
B = \left(\begin{array}{c c}
		-10 & 0 \\
		0 & -10
	\end{array}\right).
$$
It is easy to verify that 
$$D = \left(\begin{array}{c c}
		4 & 0 \\
		0 & 1
	\end{array}\right)$$
satisfies $A^TD+DA \prec0$, $(A+B)^TD+D(A+B) \prec 0$.  However, if we consider $H$ of the form \eqref{eq:H} with
$$H_{11} = \left(\begin{array}{c c}
		3 & 0 \\
		0 & 3
	\end{array}\right), \,\,
H_{12} = \left(\begin{array}{c c}
		-1 & 0 \\
		0 & -1
	\end{array}\right), \,\,
H_{22} = \left(\begin{array}{c c}
		2 & 0 \\
		0 & 2
	\end{array}\right),$$
then $H \succeq 0$ and $\textrm{diag}(H_{22}) \leq \textrm{diag}(H_{11})$.  However, it can be verified by direct computation that $AH_{11} + BH_{12}^T$ has no negative diagonal entries.  Hence, by Theorem \ref{thm:TM}, there do not exist diagonal $P \succ 0$, $Q \succ 0$ satisfying \eqref{eq:Ricc3}.
\end{example}

\section{Applications}
\label{sec:Apps}
We now present a number of applications of the general result 
Theorem \ref{thm:TM}; first considering the case when $A$ is Metzler, and $B$ is nonnegative.

\begin{theorem}
\label{thm:Metz}Let $A \in \mathbb{R}^{n \times n}$ be Metzler and $B \in \mathbb{R}^{n \times n}$ be nonnegative.  There exist diagonal matrices $P \succ 0$, $Q \succ 0$ satisfying \eqref{eq:Ricc3} if and only if $A+B$ is Hurwitz.
\end{theorem}
\textbf{Proof:} We know that if such matrices $P, Q$ exist, then there exists a diagonal $D \succ 0$ such that $(A+B)^T D + D(A+B) \prec 0$ and, in particular, that $A+B$ is Hurwitz.  

Conversely, suppose that $A+B$ is Hurwitz.  Let a non-zero
\begin{equation}
\label{eq:HUV} H = \left(\begin{array}{c c}
				U & V \\
				V^T & W
			\end{array}\right)
\end{equation}
in $Sym(2n, \mathbb{R})$ be given with $\textrm{diag}(W) \leq \textrm{diag}(U)$.  We wish to show that the matrix $AU + BV^T$ has a negative diagonal entry.  To this end, note that as $H \succeq 0$, it follows that $U \succeq 0$ and hence that
$$|u_{ji}| \leq \sqrt{u_{ii}} \sqrt{u_{jj}}$$
for $1 \leq i, j \leq n$ and of course $u_{ii} \geq 0$ for $1 \leq i \leq n$.  Furthermore, as $H \succeq 0$, we also have that 
$$|v_{ij}| \leq \sqrt{u_{ii}}\sqrt{w_{jj}} \leq \sqrt{u_{ii}}\sqrt{u_{jj}}.$$
Writing $g$ for the vector in $\mathbb{R}^n$ whose $i$th component is given by $\sqrt{u_{ii}}$, $g \neq 0$ as otherwise $H = 0$.  As $A$ is Metzler and $B$ is nonnegative, we have
\begin{eqnarray*}
(AU + BV^T)_{ii} &=& \sum_{j=1}^n a_{ij}u_{ji} + \sum_{j=1}^n b_{ij} v_{ij} \\
&\leq& \sum_{j=1}^n a_{ij} g_j g_i + \sum_{j=1}^n b_{ij} g_j g_i \\
&=& g_{i}[ (A + B)g]_i.
\end{eqnarray*}
As $A+B$ is Metzler and Hurwitz and $g \neq 0$, it follows from Proposition \ref{prop:Met1} that there exists some index $i$ such that $g_{i}[ (A + B)g]_i < 0$ and hence that 
$(AU + BV^T)_{ii}$.  It now follows from Theorem \ref{thm:TM} that there exist diagonal matrices $P \succ 0$, $Q \succ 0$ satisfying \eqref{eq:Ricc3} as claimed.

\vspace{5mm}

Our next result gives conditions under which lower triangular pairs $A, B$ are diagonally Riccati stable.

\begin{theorem}
\label{thm:LTri} Let $A \in \mathbb{R}^{n \times n}$, $B \in \mathbb{R}^{n \times n}$ be lower triangular matrices.  Then there exist diagonal matrices $P \succ 0$, $Q \succ 0$ satisfying \eqref{eq:Ricc3} if and only if $A$ is Hurwitz and $|b_{ii}| < |a_{ii}|$ for $1 \leq i \leq n$.
\end{theorem}
\textbf{Proof:} 
If there exist diagonal matrices $P \succ 0$, $Q \succ 0$ satisfying \eqref{eq:Ricc3} then $A^TD+DA \prec 0$ and hence $A$ is Hurwitz.  To show that $|b_{ii}| < |a_{ii}|$ for all $i$, we argue by contradiction. Assume that there is some $i$ with $|b_{ii}| \geq |a_{ii}|$.  First assume that $b_{ii} \geq 0$ (the case $b_{ii} \leq 0$ can  be handled identically.)  Consider the matrix $H \in Sym(2n, \mathbb{R})$ of the form \eqref{eq:HUV} where 
$$U = V = W= e_ie_i^T.$$  So $U$ and $V$ are both the diagonal matrix with all entries zero apart from the $i$th diagonal entry.  It is then easy to see that $H \succeq 0$ and that $\textrm{diag}(W) \leq \textrm{diag}(U)$.  The only non-zero diagonal entry in $AU + BV^T$ is 
\begin{eqnarray*}
(AU + BV^T)_{ii} &=& a_{ii} u_{ii} + b_{ii} v_{ii}\\
&=& a_{ii} + b_{ii} \\
&\geq& 0
\end{eqnarray*}
as $b_{ii} \geq 0$ and $|b_{ii}| \geq |a_{ii}|$.  It follows from Theorem \ref{thm:TM} that there cannot exist diagonal matrices $P \succ 0$, $Q \succ 0$ satisfying \eqref{eq:Ricc3} which is a contradiction.  For the case where $b_{ii} \leq 0$, choose $H$ with $U=W = e_ie_i^T$, $V = -e_ie_i^T$.

For the converse, assume $A$ is Hurwitz and that $|b_{ii}| < |a_{ii}|$ for $1 \leq i \leq n$.  Consider $H \succeq 0$ in $Sym(2n, \mathbb{R})$ given by \eqref{eq:HUV} with $\textrm{diag}(W) \leq \textrm{diag}(U)$.  As $H$ is non-zero, let $i$ be the lowest index such that $u_{ii} > 0$. As $U$ is positive semi-definite, it follows that $u_{ji} = 0$ for all $j < i$.  Moreover, from our assumptions on $H$, we must have 
\begin{eqnarray*}
|v_{ij}| &\leq& \sqrt{u_{ii}} \sqrt{w_{jj}}\\
&\leq& \sqrt{u_{ii}} \sqrt{u_{jj}} \\
&=& 0
\end{eqnarray*}
for all $j < i$.  Moreover, we also have $|v_{ii}| \leq |u_{ii}|$. Using these observations and the fact that $A$, $B$ are lower triangular, we see that
\begin{eqnarray*}
(AU + BV^T)_{ii} &=& \sum_{j=1}^n a_{ij} u_{ji} + \sum_{j=1}^n b_{ij}v_{ij}\\
&=& \sum_{j \leq i} a_{ij} u_{ji} + \sum_{j \leq i} b_{ij}v_{ij} \\
&=& a_{ii} u_{ii} + b_{ii} v_{ii}.
\end{eqnarray*}
The result now follows from the simple observation that $a_{ii}u_{ii} < 0$ and $|b_{ii} v_{ii}| < |a_{ii}u_{ii}|$.  

\vspace{5mm}

An identical argument where we choose the index $i$ to be equal to the largest index for which $u_{ii} > 0$ will show that a corresponding result also hold for upper triangular matrices.

\textbf{Comment:}

Note that Example \ref{ex:Ctr1} shows that even when both $A$ and $B$ are lower triangular and $A$ and $A+B$ are both Hurwitz, it is possible that the pair $A, B$ is not diagonally Riccati stable.

\section{Metzler and Nonnegative Matrices - A Closer Look}
In the previous section, we showed that if $A$ is Metzler, $B$ nonnegative and 
$A+B$ Hurwitz, there exist diagonal matrices 
$P \succ 0$, $Q \succ 0$ satisfying \eqref{eq:Ricc3}.  We have also noted 
that for any such pair $P, Q$, the matrix $P$ must be a simultaneous diagonal 
solution of the Lyapunov inequalities associated with $A$, $A+B$.  We next 
show that for any diagonal $P$ satisfying 
$A^TP+PA \prec 0$, $(A+B)^TP+P(A+B) \prec 0$, there must exist 
a \emph{diagonal} $Q$ such that $P$, $Q$ satisfy \eqref{eq:Ricc3}.  Note that in this case, as $B$ is nonnegative, from  $(A+B)^TP+P(A+B) \prec 0$
it follows that $A^TP+PA \prec 0$.
 
This strengthens the result given in \cite{MAS} and moreover, the 
method of proof used here is significantly different.  Moreover, the method of proof gives us a constructive approach to finding $P$ and $Q$.

\begin{theorem}
\label{thm:RiccDiag} Let $A \in \mathbb{R}^{n \times n}$ be Metzler and $B \in \mathbb{R}^{n \times n}$ be nonnegative and assume that $A+B$ is Hurwitz.  Let $P \succ 0$ be diagonal and such that 
\begin{equation}
\label{eq:ABdiag}
(A+B)^T P + P(A+B) \prec 0.
\end{equation}
Then there exists $Q \succ 0$, \emph{diagonal} such that $P, Q$ satisfy \eqref{eq:Ricc3}.
\end{theorem}

\textbf{Proof:}  As $P \succ 0$ is diagonal it is straightforward to verify that $(A+B)^T P + P(A+B)$ is again Metzler and, moreover as it is negative definite, it is Hurwitz.  From Proposition \ref{prop:Met1}, there exist vectors $v \gg 0$, $w \gg 0$ such that 
\begin{equation}
\label{eq:vw}
((A+B)^T P + P(A+B))v = -w.
\end{equation}  
Now choose a diagonal $Q \succ 0$ such that 
\begin{equation}
\label{eq:Q1}
Qv = B^TP v + \frac{1}{2} w.
\end{equation}
A simple calculation now shows that 
\begin{equation}
\label{eq:Q2}
\left(\begin{array}{c c}
		A^T P + PA + Q & PB \\
		B^TP & -Q
	\end{array}\right)
\left(\begin{array}{c}
		v\\
		v
	\end{array}\right) = \left(\begin{array}{c}
		-w/2\\
		-w/2
	\end{array}\right).
\end{equation}
Noting that 
$$M = \left(\begin{array}{c c}
		A^T P + PA + Q & PB \\
		B^TP & -Q
	\end{array}\right)$$ 
is Metzler, it follows from \eqref{eq:Q2} that $M$ is Hurwitz.  However as $M$ is symmetric this implies that $M \prec 0$ and hence that $P$, $Q$ satisfy \eqref{eq:Ricc3}.  This completes the proof.

\section{Generalised Lotka-Volterra Systems}
\label{sec:LV_Main}
In this section, we highlight how the results described above can be applied to a generalised class of Lotka-Volterra models occurring in population dynamics.  Specifically, we consider the following generalized Lotka--Volterra model of population dynamics in an $n$ species community:
\begin{equation}
\label{eq:GLV}
\dot x_i(t)=g_i(x_i(t))\left( c_i+\sum_{j=1}^n a_{ij} f_j(x_j(t))+
\sum_{j=1}^n 
b_{ij} f_j(x_j(t-\tau))\right), 
\end{equation} 
for $1 \leq i \leq n$.

For systems without time-delay, similar models were considered in 
\cite{RED}, while more recently a switched system version without 
time-delay was studied in \cite{ALEX1, ALEX2}.  Systems of the form \eqref{eq:GLV} with time delay and with some special forms of functions
$f_i$ were studied, for example, in \cite{FAN, CHEN}.
 
Here $x_i(t)$ is the population density of the $i$-th species; the functions
$g_i:[0, +\infty) \rightarrow [0, +\infty)$ and  $f_i:[0, +\infty) \rightarrow [0, +\infty)$
possess special properties described below and allow for a generalisation of the classical Lotka--Volterra model;
$c_i$, $a_{ij}$,  $b_{ij}$ are constant coefficients.
The coefficients $c_i$ characterise the intrinsic growth rate of the $i$-th population;
the self-interaction terms
$p_{ii}g_i(x_i) f_i(x_i)$ with $p_{ii}<0$ reflect the 
limited resources available in the environment; the terms
$p_{ij}g_i(x_i) f_j(x_j)$ for $j\neq i$ describe the influence of population $j$
on population $i$.  

Let $\mathbb{R}^n_+$ be the nonnegative cone of $\mathbb{R}^n$:
$$\mathbb{R}^n_+:=\{x \in \mathbb{R}^n \mid x \geq 0\}.$$
Also, ${\rm int} \, \mathbb{R}^n_+$ is the interior of $\mathbb{R}^n_+$.

We make use of some standard results and concepts from the Lyapunov 
stability theory of Functional Differential Equations \cite{HALE}.  
In particular, for a given  real number $\tau>0$, 
$C([-\tau, 0], {\rm int} \, \mathbb{R}^n_+)$ denotes the space of continuous 
functions 
$\phi(\theta): [-\tau,0]     \to {\rm int} \, \mathbb{R}^n_+$
with the uniform norm 
$ \|\phi\|_\tau=\max_{\theta\in [-\tau,0]}
\|\phi(\theta)\|,$  and $\|\cdot\|$ denotes the
                   Euclidean norm of a vector.
We assume that initial functions for \eqref{eq:GLV} belong to the space 
$C([-\tau, 0], {\rm int} \, \mathbb{R}^n_+)$.
For a solution $x(t)$ of \eqref{eq:GLV}, we denote by $x_t$ the element 
in $C([-\tau, 0], {\rm int} \, \mathbb{R}^n_+)$ given 
by $x_t(u) = x(t+u)$ for $u \in [-\tau, 0]$.  
For $\phi$ in $C([-\tau, 0], {\rm int} \, \mathbb{R}^n_+)$, $\phi_j$ denotes its $j$th component.

We consider functions $f_i$ and $g_i$ satisfying the following properties, which are consistent with the standard assumptions made in \cite{RED, ALEX1} and elsewhere.  Moreover, these assumptions are clearly satisfied by the classical Lotka--Volterra model \cite{HOF}.

\begin{itemize}
\item[(i)] $g_i$ and $f_i$ are continuous on $[0, +\infty)$.
\item[(ii)] For each initial condition in ${\rm int} \, \mathbb{R}^n_+$, there exists a unique, globally defined solution to \eqref{eq:GLV}. 
\item[(iii)] $g_i(0)=f_i(0)=0$,
$g_i(x_i)>0$, $f_i(x_i)>0$  for $x_i>0$.
\item[(iv)] $f_i$ is a strictly increasing function on $[0,+\infty)$, and
$f_i(x_i)\rightarrow +\infty$ as $x_i\rightarrow +\infty$.
\item[(v)] $\int_1^{+\infty}{f_i(\zeta)}/{g_i(\zeta)}\,
d\zeta=+\infty$;
\item[(vi)] $\int_0^{1} 1/{g_i(\zeta)}\,
d\zeta=+\infty.$
\end{itemize}

The assumptions on the system imply that ${\rm int} \, \mathbb{R}^n_+$  is an invariant set for \eqref{eq:GLV}.  For biological reasons, we will consider \eqref{eq:GLV} with respect to the state space ${\rm int} \, \mathbb{R}^n_+$.  The following result demonstrates the relevance of diagonal Riccati stability to Lotka--Volterra systems and provides a practical motivation for the question considered in this paper.  

Throughout the remainder of this section, $A$ and $B$ denote the matrices $A=\{a_{ij}\}_{i,j=1}^n$, $B=\{b_{ij}\}_{i,j=1}^n$.

\begin{theorem}
\label{thm:GenRes} Consider the system described by \eqref{eq:GLV}.  Assume that:
\begin{enumerate}
\item there exists an interior equilibrium point $\overline{x}$;
\item there exist diagonal matrices $P \succ 0$, $Q \succ 0$ satisfying \eqref{eq:Ricc3}. 
\end{enumerate}
Then the equilibrium point $\overline{x}$ is globally asymptotically stable 
in ${\rm int} \, \mathbb{R}^n_+$ for any value of the delay $\tau$.
\end{theorem}


\textbf{Proof:}    Denote $f(x)=(f_1(x_1),\ldots,f_n(x_n))^T$,
$c=(c_1,\ldots,c_n)^T$. 

By assumption, $\overline{x}$ is an equilibrium point in ${\rm int} \, \mathbb{R}^n_+$ so it satisfies:
$$
c+(A+B)f(x)=0. 
$$

Using this, the system \eqref{eq:GLV} can be rewritten as follows
\begin{eqnarray}
\nonumber 
\dot x_i(t) &=& g_i(x_i(t))\left( \sum_{j=1}^n a_{ij} (f_j(x_j(t))-f_j(\bar x_j)) \right.\\
\label{eq:sys2}
&+& \left. \sum_{j=1}^n 
b_{ij} (f_j(x_j(t-\tau))-f_j(\bar x_j))\right), \quad i=1,\ldots,n.
\end{eqnarray}
 
We show that for the system \eqref{eq:sys2}, there exists a Lyapunov-Krasovskii functional of the form

\begin{equation}
\label{eq:LK1}
V(\phi)=\sum_{i=1}^n p_i  \int_{\bar x_i}^{\phi(0)} 
\frac{f_i(\zeta)-f_i(\bar x_i)}{g_i(\zeta)}d\zeta 
 +\sum_{j=1}^n \mu_{j} 
\int_{-\tau}^0 (f_j (\phi_j(u))-f_j(\bar x_j))^2du, 
\end{equation}                         
where $p_i$ and $\mu_{j}$ are positive coefficients.  

In fact, we know by assumption that there exist diagonal matrices $P \succ 0$, $Q \succ 0$ such that
$$M = \left(\begin{array}{c c}
		A^T P + PA + Q & PB \\
		B^TP & -Q
	\end{array}\right) \prec - 2 \beta I
$$
for some $\beta >0$.
It follows readily that for $x, y \in \mathbb{R}^n$, 
\begin{equation}
\label{eq:xy} 
x^TPA x + x^TPBy + 1/2 (x^TQx - y^TQy) < -\beta (x^Tx + y^T y).
\end{equation}
Choosing $p_1, \ldots, p_n$ and $\mu_1, \ldots, \mu_n$ so that $P = \textrm{diag}(p_1, \ldots, p_n)$, $Q = \rm{diag}(2\mu_1, \ldots, 2\mu_n)$, we claim that $V$ given by \eqref{eq:LK1} defines a Lyapunov-Krasovskii functional for the system \eqref{eq:sys2}.  

First of all, it follows from the assumed properties of the functions $f_i$ and $g_i$ that
$$
\chi_1(\|\varphi(0)-\bar x\|)\leq V(\varphi)
\leq \chi_2(\|\varphi-\bar x\|_\tau),
$$
where $\chi_1(\zeta)$ and $\chi_2(\zeta)$ are Khan functions \cite{ROUCHE},
and 
$V(\varphi)\to +\infty$ as $\varphi(0)$ tends to the boundary of 
                 ${\rm int} \, \mathbb{R}^n_+$. 

Differentiating the functional \eqref{eq:LK1} with respect to the system \eqref{eq:sys2}, we obtain
$$
\dot V(x_t)=\sum_{i=1}^n p_i   
(f_i(x_i(t))-f_i(\bar x_i)) 
\sum_{j=1}^n a_{ij} (f_j(x_j(t))-f_j(\bar x_j))
$$
  $$
+\sum_{i=1}^n p_i   
(f_i(x_i(t))-f_i(\bar x_i))     \sum_{j=1}^n 
b_{ij} (f_j(x_j(t-\tau))-f_j(\bar x_j))
 $$
$$
 +\sum_{j=1}^n \mu_{j} 
(f_j (x_j(t))-f_j(\bar x_j))^2   -
\sum_{j=1}^n \mu_{j} 
(f_j(x_j(t-\tau))-f_j(\bar x_j))^2.
$$     

Writing this more compactly, we have 
\begin{eqnarray*}
\dot V(x_t) = (f(x(t)) - f(\bar{x}))^T P A(f(x(t)) - f(\bar{x}))  
\\ +(f(x(t)) - f(\bar{x}))^T P B(f(x(t - \tau)) - f(\bar{x}))\\
+ \frac12(f(x(t)) - f(\bar{x}))^T Q(f(x(t)) - f(\bar{x})) \\
+ \frac12(f(x(t-\tau)) - f(\bar{x}))^T Q(f(x(t-\tau)) - f(\bar{x})).
\end{eqnarray*}
It follows immediately from \eqref{eq:xy} that
$$
\dot V(x_t)\leq -\beta
\sum_{i=1}^n   
\left((f_i(x_i(t))-f_i(\bar x_i))^2 
+(f_i (x_i(t-\tau))-f_i(\bar x_i))^2\right).
$$      This completes the proof.

\vspace{5mm}

\textbf{Comment:} The previous result shows that if a generalised Lotka-Volterra system has an interior equilibrium point and if there exist diagonal solutions $P \succ 0$, $Q \succ 0$ to  \eqref{eq:Ricc3} then the equilibrium is globally asymptotically stable independent of delay.  Moreover, the stability is established using a simple \emph{diagonal} Lyapunov-Krasovskii functional.  We next describe the implication of this result for so-called \emph{mutualistic systems} in which the matrix $A$ is Metzler and $B$ is nonnegative.

\begin{corollary}
\label{cor:Mut} Consider the system \eqref{eq:GLV}.  Assume that $A$ 
is Metzler, $B$ is nonnegative and that $c \gg 0$.  Furthermore, assume 
that $A+B$ is Hurwitz.  Then there exists an equilibrium point $\bar{x}$ in 
${\rm int} \, \mathbb{R}^n_+$ which is globally asymptotically stable for all values of the delay $\tau$. 
\end{corollary}

\textbf{Proof:} As $A+B$ is Metzler and Hurwitz and $c \gg 0$, it follows that 
$$-(A+B)^{-1} c \gg 0.$$
It now follows immediately from property (iv) of the functions 
$f_i$, $1\leq i \leq n$, that there exists a unique $\bar{x} \in 
{\rm int} \, \mathbb{R}^n_+$ satisfying 
$$f(\bar{x}) = -(A+B)^{-1} c,$$
which implies that $\bar{x}$ is an equilibrium of \eqref{eq:GLV}.  Moreover, as $A+B$ is Metzler and Hurwitz, Theorem \ref{thm:Metz} implies that there exist diagonal matrices $P \succ 0$, $Q \succ 0$ satisfying \eqref{eq:Ricc3}.  The result now follows from Theorem \ref{thm:GenRes}.

\vspace{5mm}

To finish off this section, we note that the existence of diagonal $P$, $Q$ satisfying the Riccati inequality also implies uniform boundedness of the generalised system \eqref{eq:GLV} in the case where there is not an interior equilibrium point.  

The system \eqref{eq:GLV} is uniformly ultimately bounded with the ultimate  bound $R>0$ if,
for any  $M>0$, there exists $T=T(M)> 0$, such that
  $\|{\bf x}(t, \phi,t_0)\|\leq R$ for all
 $t_0\geq 0$,  $\|\phi\|_\tau<M$,  $t\geq t_0+T$.

Choose an interior point, not necessarily an equilibrium, $\overline{x} 
\in {\rm int} \, \mathbb{R}^n_+$. Consider the Lyapunov-Krasovskii
 functional \eqref{eq:LK1}.

As in the proof of Theorem \ref{thm:GenRes}, we obtain that
$V(\phi)\to +\infty$ as $\phi(0)$ tends to the boundary of 
                 ${\rm int} \, \mathbb{R}^n_+$. 

A relatively straightforward calculation shows that for the derivative of the functional \eqref{eq:LK1} 
with respect to the system \eqref{eq:GLV} the estimate 
$$
\dot V(x_t)\leq -\beta
\sum_{i=1}^n   
\left((f_i(x_i(t))-f_i(\bar x_i))^2 
+(f_i (x_i(t-\tau))-f_i(\bar x_i))^2\right)
$$
$$
+\eta\sum_{i=1}^n   
\left|f_i(x_i(t))-f_i(\bar x_i)\right|  \eqno(19)
$$      
holds, where $\beta$ and $\eta$ are positive constants.

This implies that there exists a number $L>0$ such that $\dot V(x_t)<0$ for $\|x(t)\|>L$,
$x(t)\in \mathbb{R}^n_+$, which in turn implies that the system \eqref{eq:GLV} is uniformly ultimately bounded in $\mathbb{R}^n_+$.

\section{Concluding Remarks}
\label{sec:conc}
We have considered the general question of when a pair of matrices $A, B$ is diagonally Riccati stable.  Building on the result given in \cite{MAS}, we have derived a characterisation of diagonally Riccati stable pairs the echoes the classical result of Barker, Berman and Plemmons for the Lyapunov equation.  Applications to Metzler and Triangular matrices are presented and the relevance of diagonal Riccati stability to generalised Lotka-Volterra systems in population dynamics has been highlighted. 

\section*{Acknowledgements}

The first named author's  research was supported
by the Saint Petersburg State University (project no. 9.38.674.2013), and by the Russian Foundation of Basic Research (grant
Nos. 13-01-00347-a and 13-01-00376-a).

\end{document}